\documentclass[onecolumn]{autart}    % Enable this line and disable the
\usepackage{graphicx}          % Include this line if your
\begin{document}
\input{amssym}
\begin{frontmatter}

\title{Symmetry Analysis of Telegraph Equation}
%\thanksref{footnoteinfo}} % Title, preferably not more

%\thanks[footnoteinfo]{This paper was not presented at any IFAC
%meeting. Corresponding author M.~T.~Cicero. Tel. +XXXIX-VI-mmmxxi.
%Fax +XXXIX-VI-mmmxxv.}

\author[]{Mehdi Nadjafikhah}\ead{m\_nadjafikhah@iust.ac.ir},    % Add the
\author[MN]{Seyed Reza Hejazi}\ead{reza\_hejazi@iust.ac.ir},               % e-mail address
%\author[Baiae]{Publius Maro Vergilius}\ead{vergilius@culture.ir}  % (ead) as shown

\address[MN]{School of Mathematics, Iran University of Science and Technology, Narmak-16, Tehran, I.R.Iran}  % Please supply
%\address[Iran]{Iran University of Science and Technology, Narmak, Tehran}             % full addresses
%\address[Baiae]{The White House, Baiae}        % here.

%
% M.S.C. 2000: 49J30, 49K30, 55S45.
%
                             % keyword list or with the
                                          % help of the Automatica
                                          % keyword wizard

\begin{abstract}                          % Abstract of not more than 200 words.
Lie symmetry group method is applied to study the Telegraph
equation. The symmetry group and its optimal system are given, and
group invariant solutions associated to the symmetries are
obtained. Finally the structure of the Lie algebra symmetries is
determined.
\end{abstract}

\end{frontmatter}
%-------------------------------------------------------------------------------------------------
%-------------------------------------------------------------------------------------------------
\section{Introduction}
The telegrapher's equations (or just telegraph equations) are a
pair of linear differential equations which describe the voltage
and current on an electrical transmission line with distance and
time. The equations come from \textit{Oliver Heaviside} who
developed the transmission line model. Oliver Heaviside (May 18,
1850 – February 3, 1925) was a self-taught English electrical
engineer, mathematician, and physicist who adapted complex numbers
to the study of electrical circuits, invented mathematical
techniques to the solution of differential equations (later found
to be equivalent to Laplace transforms), reformulated Maxwell's
field equations in terms of electric and magnetic forces and
energy flux, and independently co-formulated vector analysis.
Although at odds with the scientific establishment for most of his
life, Heaviside changed the face of mathematics and science for
years to come the theory applies to high-frequency transmission
lines (such as telegraph wires and radio frequency conductors) but
is also important for designing high-voltage energy transmission
lines. The model demonstrates that the electromagnetic waves can
be reflected on the wire, and that wave patterns can appear along
the line. The telegrapher's equations can be understood as a
simplified case of Maxwell's equations. In a more practical
approach, one assumes that the conductors are composed of an
infinite series of two-port elementary components, each
representing an infinitesimally short segment of the transmission
line.
\section{Lie Symmetries of the Equation}
A PDE with $p-$independent and $q-$dependent variables has a Lie
point transformations
\begin{eqnarray*}
\widetilde{x}_i=x_i+\varepsilon\xi_i(x,u)+{\mathcal
O}(\varepsilon^2),\qquad
\widetilde{u}_{\alpha}=u_\alpha+\varepsilon\varphi_\alpha(x,u)+{\mathcal
O}(\varepsilon^2)
\end{eqnarray*}
where
$\displaystyle{\xi_i=\frac{\partial\widetilde{x}_i}{\partial\varepsilon}\Big|_{\varepsilon=0}}$
for $i=1,...,p$ and
$\displaystyle{\varphi_\alpha=\frac{\partial\widetilde{u}_\alpha}{\partial\varepsilon}\Big|_{\varepsilon=0}}$
for $\alpha=1,...,q$. The action of the Lie group can be
considered by its associated infinitesimal generator
\begin{eqnarray}\label{eq:18}
\textbf{v}=\sum_{i=1}^p\xi_i(x,u)\frac{\partial}{\partial{x_i}}+\sum_{\alpha=1}^q\varphi_\alpha(x,u)\frac{\partial}{\partial{u_\alpha}}
\end{eqnarray}
on the total space of PDE (the space containing independent and
dependent variables). Furthermore, the characteristic of the
vector field (\ref{eq:18}) is given by
\begin{eqnarray*}
Q^\alpha(x,u^{(1)})=\varphi_\alpha(x,u)-\sum_{i=1}^p\xi_i(x,u)\frac{\partial
u^\alpha}{\partial x_i},
\end{eqnarray*}
and its $n-$th prolongation is determined by
\begin{eqnarray*}
\textbf{v}^{(n)}=\sum_{i=1}^p\xi_i(x,u)\frac{\partial}{\partial
x_i}+\sum_{\alpha=1}^q\sum_{\sharp
J=j=0}^n\varphi^J_\alpha(x,u^{(j)})\frac{\partial}{\partial
u^\alpha_J},
\end{eqnarray*}
where
$\varphi^J_\alpha=D_JQ^\alpha+\sum_{i=1}^p\xi_iu^\alpha_{J,i}$.
($D_J$ is the total derivative operator describes in
(\ref{eq:19})).

The aim is to analysis the point symmetry structure of the
Telegraph equation, which is
\begin{equation}\label{eq:1}
u_{tt}+ku_t=a^2\Big[\frac{1}{r}\Big(r\frac{\partial u}{\partial
r}\Big)_r+\frac{1}{r^2}u_{xx}+u_{yy}\Big],
\end{equation}
where $u$ is a smooth function of $\displaystyle{(r,x,y,t)}$.

Let us consider a one-parameter Lie group of infinitesimal
transformations $(x,t,u)$ given by
\begin{eqnarray*}\begin{array}{lllll}
\widetilde{r}=r+\varepsilon\xi_1(r,x,y,t,u)+{\mathcal
O}(\varepsilon^2),&&
\widetilde{x}=x+\varepsilon\xi_2(r,x,y,t,u)+{\mathcal
O}(\varepsilon^2),&&
\widetilde{y}=y+\varepsilon\xi_3(r,x,y,t,u)+{\mathcal
O}(\varepsilon^2),\\
\widetilde{t}=t+\varepsilon\xi_4(r,x,y,t,u)+{\mathcal
O}(\varepsilon^2),&&
\widetilde{u}=u+\varepsilon\eta(r,x,y,t,u)+{\mathcal
O}(\varepsilon^2),\end{array}
\end{eqnarray*}
where $\varepsilon$ is the group parameter. Then one requires that
this transformations leaves invariant the set of solutions of the
Eq. (\ref{eq:1}). This yields to the linear system of equations
for the infinitesimals $\xi_1(r,x,y,t,u)$, $\xi_2(r,x,y,t,u)$,
$\xi_3(r,x,y,t,u)$, $\xi_4(r,x,y,t,u)$ and $\eta(r,x,y,t,u)$. The
Lie algebra of infinitesimal symmetries is the set of vector
fields in the form of
$\textbf{v}=\xi_1\partial_r+\xi_2\partial_x+\xi_3\partial_y+\xi_4\partial_t+\eta\partial_u$.
This vector field has the second prolongation
\begin{eqnarray*}
\textbf{v}^{(2)}=\textbf{v}+\varphi^r\partial_{r}+\varphi^x\partial_{x}+\varphi^y\partial_{y}+\varphi^t\partial_{t}+\varphi^{rr}\partial_{u_{rr}}
+\varphi^{rx}\partial_{u_{rx}}+\cdots+\varphi^{yy}\partial_{u_{yy}}+\varphi^{yt}\partial_{u_{yt}}+\varphi^{tt}\partial_{tt}
\end{eqnarray*}
with the coefficients
\begin{eqnarray*}\begin{array}{lll}
\varphi^r =D_rQ+\xi_1u_{rr}+\xi_2u_{rx}+\xi_3u_{ry}+\xi_4u_{rt},&&
\varphi^x =D_xQ+\xi_1u_{rx}+\xi_2u_{xx}+\xi_3u_{xy}+\xi_4u_{xt},\\
\varphi^y =D_yQ+\xi_1u_{ry}+\xi_2u_{xy}+\xi_3u_{yy}+\xi_4u_{yt},&&
\varphi^t =D_xQ+\xi_1u_{rt}+\xi_2u_{xt}+\xi_3u_{xt}+\xi_4u_{tt},\\
\varphi^{rr}=D^2_rQ+\xi_1u_{rrr}+\xi_2u_{rrx}+\xi_3u_{rry}+\xi_4u_{rrt},&&
\varphi^{rx}=D_rD_xQ+\xi_1u_{rxr}+\xi_2u_{rxx}+\xi_3u_{rxy}+\xi_4u_{rxt},\\
\varphi^{ry}=D_rD_yQ+\xi_1u_{ryr}+\xi_2u_{rxy}+\xi_3u_{ryy}+\xi_4u_{ryt}&&
\varphi^{rt}=D_rD_tQ+\xi_1u_{rtr}+\xi_2u_{rxt}+\xi_3u_{ryt}+\xi_4u_{rtt},\\
\varphi^{xx}=D^2_xQ+\xi_1u_{xxr}+\xi_2u_{xxx}+\xi_3u_{xxy}+\xi_4u_{xxt},&&
\varphi^{xy}=D_xD_yQ+\xi_1u_{xyr}+\xi_2u_{xxy}+\xi_3u_{xyy}+\xi_4u_{xyt},\\
\varphi^{xt}=D_xD_tQ+\xi_1u_{xrt}+\xi_2u_{xxt}+\xi_3u_{xyt}+\xi_4u_{xtt},&&
\varphi^{yy}=D^2_yQ+\xi_1u_{ryy}+\xi_2u_{xyy}+\xi_3u_{yyy}+\xi_4u_{yyt},\\
\varphi^{yt}=D_yD_tQ+\xi_1u_{ryt}+\xi_2u_{xyt}+\xi_3u_{yyt}+\xi_4u_{ytt},&&
\varphi^{tt}=D^2_tQ+\xi_1u_{rtt}+\xi_2u_{xtt}+\xi_3u_{ytt}+\xi_4u_{ttt},
\end{array}
\end{eqnarray*}
where the operators $D_r,D_x,D_y$ and $D_t$ denote the total
derivative with respect to $r,x,y$ and $t$:
\begin{eqnarray}\label{eq:19}\begin{array}{lll}
D_r=\partial_r+u_r\partial_u+u_{rr}\partial_{u_r}+u_{rx}\partial_{u_x}+\cdots,&&
D_x=\partial_x+u_x\partial_u+u_{xx}\partial_{u_x}+u_{rx}\partial_{u_r}+\cdots,\\
D_y=\partial_y+u_y\partial_u+u_{yy}\partial_{u_y}+u_{ry}\partial_{u_r}+\cdots,&&
D_t=\partial_t+u_t\partial_u+u_{tt}\partial_{u_t}+u_{rt}\partial_{u_r}+\cdots,\end{array}
\end{eqnarray}
Using the invariance condition, i.e., applying the second
prolongation $\textbf{v}^{(2)}$ to Eq. (\ref{eq:1}), the following
system of 27 determining equations yields:
\begin{eqnarray*}
\begin{array}{lclclclc}
{\xi_2}_u=0,&&{\xi_2}_{yy}=0,&&{\xi_3}_y=0,&&\hspace{-5cm}{\xi_3}_u=0,\\
{\xi_4}_t=0,&&{\xi_4}_u=0,&&{\xi_4}_{rr}=0,&&\hspace{-5cm}{\xi_4}_{xy}=0,\\
{\xi_4}_{yy}=0,&&{\xi_4}_{ry}=0,&&\eta_{tu}=0,&&\hspace{-5cm}\eta_{uu}=0\\
k{\xi_4}_y+2\eta_{ru}=0,&&\xi_1+r{\xi_2}_x=0,&&{\xi_2}_x+r{\xi_2}_{rx}=0,&&\hspace{-5cm}{\xi_2}_y+r{\xi_2}_{ry}=0,\\
{\xi_2}_{xx}-r{\xi_2}_r=0,&&k{\xi_4}_y+2\eta_{yu}=0,&&2{\xi_2}_r+r{\xi_2}_{rr}=0,&&\hspace{-5cm}{\xi_3}_r-r{\xi_2}_{xy}=0,\\
{\xi_3}_x+r^2{\xi_2}_y=0,&&{\xi_3}_t-a^2{\xi_4}_y=0,&&{\xi_4}_x-r{\xi_4}_{rx}=0,&&\hspace{-5cm}{\xi_4}_{xx}+r{\xi_4}_r=0,\\
r^2{\xi_2}_t-a^2{\xi_4}_x=0,&&k{\xi_4}_x+2\eta_{ux}=0,&&a^2r^2\eta_{rr}-kr^2\eta_t+a^2r\eta_r+a^2r^2\eta_{yy}-r^2\eta_{tt}+a^2\eta_{xx}=0.
\end{array}
\end{eqnarray*}
The solution of the above system gives the following coefficients
of the vector field $\textbf{v}$:
\begin{eqnarray*}
\xi_1&=&c_6\sin x-c_7\cos x-c_8y\cos x-c_9y\sin x+2c_{10}a^2t\sin x-2c_{11}a^2t\cos x,\\
\xi_2&=&c_1+c_6r^{-1}\cos x+c_7r^{-1}\sin x+c_8yr^{-1}\sin x-c_9yr^{-1}\sin x+2c_{10}a^2tr^{-1}\cos x-2c_{11}a^2tr^{-1}\sin x,\\
\xi_3&=&c_2+2c_5a^2t+c_8r\cos x+c_9r\sin x,\\
\xi_4&=&c_3+2c_5at+2c_{10}r\sin x-2c_{11}r\cos x,\qquad
\eta=c_4u-c_5kyu-c_{10}kru\sin x+c_{11}kru\cos x,
\end{eqnarray*}
where $c_1,...,c_{11}$ are arbitrary constants, thus the Lie
algebra ${\goth g}$ of the telegraph equation is spanned by the
seven vector fields
\begin{eqnarray*}
\begin{array}{lclc}
\textbf{v}_1=\partial_x,\qquad\textbf{v}_2=\partial_y,&&\hspace{-3cm}
\textbf{v}_3=\partial_t,\qquad\textbf{v}_4=u\partial_u,\\
\textbf{v}_5=2a^2t\partial_y+2y\partial_t-kyu\partial_u,&&\hspace{-3cm}\textbf{v}_6=\sin x\partial_r+r^{-1}\cos x\partial_x,\\
\textbf{v}_7=-\cos x\partial_r+r^{-1}\sin
x\partial_x,&&\hspace{-3cm} \textbf{v}_8=-y\cos
x\partial_r+r^{-1}y\sin x\partial_x+r\cos
x\partial_y,\\\textbf{v}_9=-y\sin x\partial_r-r^{-1}y\cos
x\partial_x+r\sin
x\partial_y,&&\hspace{-3cm}\textbf{v}_{10}=2a^2t\sin
x\partial_r+2a^2tr^{-1}\cos x\partial_x+2r\sin x\partial_t-kru\sin
x\partial_u,\\\textbf{v}_{11}=-2a^2t\cos
x\partial_r+2a^2tr^{-1}\sin x\partial_x-2r\cos x\partial_t+kru\cos
x\partial_u,
\end{array}
\end{eqnarray*}
which $\textbf{v}_1,\textbf{v}_2$ and $\textbf{v}_3$ are
translation on $x,t$ and $u$, $\textbf{v}_4$ is rotation on $u$
and $x$ and $\textbf{v}_7$ is scaling on $x,t$ and $u$. The
commutation relations between these vector fields is given by the
(\ref{table:1}), where entry in row $i$ and column $j$
representing $[\textbf{v}_i,\textbf{v}_j]$.
\begin{table}
\caption{Commutation relations of $\goth g$ }\label{table:1}
\vspace{-0.3cm}\begin{eqnarray*}\hspace{-0.75cm}\begin{array}{cccccccccccc}
\hline
  [\,,\,]        &\hspace{1cm}\textbf{v}_1       &\hspace{0.5cm}\textbf{v}_2                    &\hspace{0.5cm}\textbf{v}_3  &\hspace{0.5cm}\textbf{v}_4 &\hspace{0.5cm}\textbf{v}_5                   &\hspace{0.5cm}\textbf{v}_6                     &\hspace{0.5cm}\textbf{v}_7                    &\hspace{0.5cm}\textbf{v}_8        &\hspace{0.5cm}\textbf{v}_9       &\hspace{0.5cm}\textbf{v}_{10}    &\hspace{0.5cm}\textbf{v}_{11}  \\ \hline
  \textbf{v}_1   &\hspace{1cm} 0                 &\hspace{0.5cm} 0                              &\hspace{0.5cm}0             &\hspace{0.5cm}0            &\hspace{0.5cm}0                              &\hspace{0.3cm}-\textbf{v}_7                    &\hspace{0.5cm}\textbf{v}_6                    &\hspace{0.3cm}-\textbf{v}_9       &\hspace{0.5cm}\textbf{v}_8       &\hspace{0.3cm}-\textbf{v}_{11}   &\hspace{0.5cm}\textbf{v}_{10}     \\
  \textbf{v}_2   &\hspace{1cm} 0                 &\hspace{0.5cm} 0                              &\hspace{0.5cm}0             &\hspace{0.5cm}0            &\hspace{0.5cm}2a^2\textbf{v}_3-k\textbf{v}_4 &\hspace{0.5cm}0                                &\hspace{0.5cm}0                               &\hspace{0.5cm}\textbf{v}_7        &\hspace{0.3cm}-\textbf{v}_6      &\hspace{0.5cm}0                  &\hspace{0.5cm}0\\
  \textbf{v}_3   &\hspace{1cm} 0                 &\hspace{0.5cm} 0                              &\hspace{0.5cm}0             &\hspace{0.5cm}0            &\hspace{0.5cm}\textbf{v}_2                   &\hspace{0.5cm}0                                &\hspace{0.5cm}0                               &\hspace{0.5cm}0                   &\hspace{0.5cm}0                  &\hspace{0.5cm}\textbf{v}_6       &\hspace{0.5cm}\textbf{v}_7\\
  \textbf{v}_4   &\hspace{1cm} 0                 &\hspace{0.5cm} 0                              &\hspace{0.5cm}0             &\hspace{0.5cm}0            &\hspace{0.5cm}0                              &\hspace{0.5cm}0                                &\hspace{0.5cm}0                               &\hspace{0.5cm}0                   &\hspace{0.5cm}0                  &\hspace{0.5cm}0                  &\hspace{0.5cm}0\\
  \textbf{v}_5   &\hspace{1cm} 0                 &\hspace{-0.1cm} -2a^2\textbf{v}_3+\textbf{v}_4&\hspace{0.3cm}-\textbf{v}_2 &\hspace{0.5cm}0            &\hspace{0.5cm}0                              &\hspace{0.5cm}0                                &\hspace{0.5cm}0                               &\hspace{0.5cm}\textbf{v}_{11}     &\hspace{0.3cm}-\textbf{v}_{10}   &\hspace{0.3cm}-\textbf{v}_9      &\hspace{0.5cm}\textbf{v}_8\\
  \textbf{v}_6   &\hspace{1cm} \textbf{v}_7      &\hspace{0.5cm} 0                              &\hspace{0.5cm}0             &\hspace{0.5cm}0            &\hspace{0.5cm}0                              &\hspace{0.5cm}0                                &\hspace{0.5cm}0                               &\hspace{0.5cm}\textbf{v}_7        &\hspace{0.3cm}-\textbf{v}_6      &\hspace{0.5cm}0                  &\hspace{0.5cm}0\\
  \textbf{v}_7   &\hspace{0.9cm} -\textbf{v}_6   &\hspace{0.5cm} 0                              &\hspace{0.5cm}0             &\hspace{0.5cm}0            &\hspace{0.5cm}0                              &\hspace{0.5cm}0                                &\hspace{0.5cm}0                               &\hspace{0.3cm}-\textbf{v}_2       &\hspace{0.5cm}0                  &\hspace{0.5cm}0                  &\hspace{-0.2cm}2a^2\textbf{v}_3-k\textbf{v}_4\\
  \textbf{v}_8   &\hspace{1cm} \textbf{v}_9      &\hspace{0.3cm} -\textbf{v}_7                  &\hspace{0.5cm}0             &\hspace{0.5cm}0            &\hspace{0.3cm}-\textbf{v}_{11}               &\hspace{0.5cm}0                                &\hspace{0.3cm}-\textbf{v}_2                   &\hspace{0.5cm}0                   &\hspace{0.3cm}-\textbf{v}_1      &\hspace{0.5cm}0                  &\hspace{0.5cm}\textbf{v}_5\\
  \textbf{v}_9   &\hspace{0.9cm} -\textbf{v}_8   &\hspace{0.5cm} \textbf{v}_6                   &\hspace{0.5cm}0             &\hspace{0.5cm}0            &\hspace{0.5cm}\textbf{v}_{10}                &\hspace{0.3cm}-\textbf{v}_2                    &\hspace{0.5cm}0                               &\hspace{0.5cm}\textbf{v}_1        &\hspace{0.5cm}0                  &\hspace{0.3cm}-\textbf{v}_5      &\hspace{0.5cm}0\\
  \textbf{v}_{10}&\hspace{1cm} \textbf{v}_{11}   &\hspace{0.5cm} 0                              &\hspace{0.3cm}-\textbf{v}_6 &\hspace{0.5cm}0            &\hspace{0.3cm}-\textbf{v}_9                  &\hspace{-0.4cm}-2a^2\textbf{v}_3+k\textbf{v}_4 &\hspace{0.5cm}0                               &\hspace{0.5cm}0                   &\hspace{0.5cm}\textbf{v}_5       &\hspace{0.5cm}0                  &\hspace{0.5cm}a^2\textbf{v}_1\\
  \textbf{v}_{11}&\hspace{0.9cm} -\textbf{v}_{10}&\hspace{0.5cm} 0                              &\hspace{0.3cm}-\textbf{v}_7 &\hspace{0.5cm}0            &\hspace{0.3cm}-\textbf{v}_8                  &\hspace{0.5cm}0                                &\hspace{-0.3cm}-2a^2\textbf{v}_3+k\textbf{v}_4&\hspace{0.3cm}-\textbf{v}_5       &\hspace{0.5cm}0                  &\hspace{0.3cm}-a^2\textbf{v}_1   &\hspace{0.5cm}0\\
  \hline\end{array}\end{eqnarray*}\end{table}
The one-parameter groups $G_i$ generated by the base of $\goth g$
are given in the following table.
\begin{eqnarray*}
g_1   &:&(r,x,y,t,u)\longmapsto(r,x+s,y,t,u),\qquad\qquad
g_2   :(r,x,y,t,u)\longmapsto(r,x,y+s,t,u),\\
g_3   &:&(r,x,y,t,u)\longmapsto(r,x,y,t+s,u),\qquad\qquad
g_4   :(r,x,y,t,u)\longmapsto(r,x,y,t,ue^s),\\
g_5   &:&(r,x,y,t,u)\longmapsto\Big(r,x,y+st,\frac{1}{a^2}sy+t,-\frac{1}{2a^2}skyu+u\Big),\\
g_6   &:&(r,x,y,t,u)\longmapsto\Big(s\sin x+r,\frac{s}{r}\cos x+x,y,t,u\Big),\\
g_7   &:&(r,x,y,t,u)\longmapsto\Big(-s\cos x+r,\frac{s}{r}\sin x+x,y,t,u\Big),\\
g_8   &:&(r,x,y,t,u)\longmapsto\Big(-sy\cos x+r,\frac{s}{r}y\sin x+x,sr\cos x+y,t,u\Big),\\
g_9   &:&(r,x,y,t,u)\longmapsto\Big(-sy\sin x+r,-\frac{s}{r}y\cos x+x,sr\sin x+y,t,u\Big),\\
g_{10}&:&(r,x,y,t,u)\longmapsto\Big(st\sin x+r,\frac{s}{r}t\cos x+x,y,\frac{s}{a^2}r\sin x+t-\frac{s}{2a^2}kru\sin x+u\Big),\\
g_{11}&:&(r,x,y,t,u)\longmapsto\Big(-st\cos x+r,\frac{s}{r}t\sin
x+x,y,-\frac{s}{a^2}r\cos x+t,\frac{s}{2a^2}kru\cos x+u\Big).
\end{eqnarray*}
Since each group $G_i$ is a symmetry group and if $u=f(r,x,y,t)$
is a solution of the Telegraph equation, so are the functions
\begin{eqnarray*}
u^1&=&U(r,x+\varepsilon,y,t),\qquad
u^2=U(r,x,y+\varepsilon,t),\qquad
u^3=U(r,x,y,t+\varepsilon),\qquad
u^4=e^{-\varepsilon}U(r,x,y,t),\\
u^5   &=&(2a^2+\varepsilon ky)U\Big(r,x,y+\varepsilon
t,\frac{1}{a^2}\varepsilon y+t\Big),\qquad
u^6   =U\Big(\varepsilon\sin x+r,\frac{\varepsilon}{r}\cos x+x,y,t\Big),\\
u^7   &=&U\Big(-\varepsilon\cos x+r,\frac{\varepsilon}{r}\sin
x+x,y,t\Big),\qquad\quad\;
u^8=U\Big(-\varepsilon y\cos x+r,\frac{\varepsilon}{r}y\sin x+x,\varepsilon r\cos x+y,t\Big),\\
u^9   &=&U\Big(-\varepsilon y\sin x+r,-\frac{\varepsilon}{r}y\cos x+x,\varepsilon r\sin x+y,t\Big),\\
u^{10}&=&(2a^2+\varepsilon kr\sin x)U\Big(\varepsilon t\sin x+r,\frac{\varepsilon}{r}t\cos x+x,y,\frac{\varepsilon}{a^2}r\sin x+t\Big),\\
u^{11}&=&(2a^2-\varepsilon kr\cos x)U\Big(-\varepsilon t\cos
x+r,\frac{\varepsilon}{r}t\sin x+x,y,-\frac{\varepsilon}{a^2}r\cos
x+t\Big).
\end{eqnarray*}
where $\varepsilon$ is a real number. Here we can find the general
group of the symmetries by considering a general linear
combination $c_1\textbf{v}_1+\cdots+c_1\textbf{v}_{11}$ of the
given vector fields. In particular if $g$ is the action of the
symmetry group near the identity, it can be represented in the
form
$g=\exp(\varepsilon_{11}\textbf{v}_{11})\cdots\exp(\varepsilon_1\textbf{v}_1)$.
\section{Optimal system of Telegraph equation}
$~~~~~$As is well known, the theoretical Lie group method plays an
important role in finding exact solutions and performing symmetry
reductions of differential equations. Since any linear combination
of infinitesimal generators is also an infinitesimal generator,
there are always infinitely many different symmetry subgroups for
the differential equation. So, a mean of determining which
subgroups would give essentially different types of solutions is
necessary and significant for a complete understanding of the
invariant solutions. As any transformation in the full symmetry
group maps a solution to another solution, it is sufficient to
find invariant solutions which are not related by transformations
in the full symmetry group, this has led to the concept of an
optimal system \cite{[6]}. The problem of finding an optimal
system of subgroups is equivalent to that of finding an optimal
system of subalgebras. For one-dimensional subalgebras, this
classification problem is essentially the same as the problem of
classifying the orbits of the adjoint representation. This problem
is attacked by the naive approach of taking a general element in
the Lie algebra and subjecting it to various adjoint
transformations so as to simplify it as much as possible. The idea
of using the adjoint representation for classifying
group-invariant solutions is due to \cite{[1-1],[2-1],[5],[6]}.

The adjoint action is given by the Lie series
\begin{eqnarray}\label{eq:9}
\mbox{Ad}(\exp(\varepsilon\textbf{v}_i)\textbf{v}_j)=\textbf{v}_j-\varepsilon[\textbf{v}_i,\textbf{v}_j]+\frac{\varepsilon^2}{2}[\textbf{v}_i,[\textbf{v}_i,\textbf{v}_j]]-\cdots,
\end{eqnarray}
where $[\textbf{v}_i,\textbf{v}_j]$ is the commutator for the Lie
algebra, $\varepsilon$ is a parameter, and $i,j=1,\cdots,11$. Let
$F^{\varepsilon}_i:{\goth g}\rightarrow{\goth g}$ defined by
$\textbf{v}\mapsto\mbox{Ad}(\exp(\varepsilon\textbf{v}_i)\textbf{v})$
is a linear map, for $i=1,\cdots,11$. The matrices
$M^\varepsilon_i$ of $F^\varepsilon_i$, $i=1,\cdots,11$, with
respect to basis $\{\textbf{v}_1,\cdots,\textbf{v}_{11}\}$ are
\begin{eqnarray*}
&\displaystyle
M^\varepsilon_1=\tiny\left(\begin{array}{ccccccccccc}
1&0&0&0&0&0&0&0&0&0&0\\0&1&0&0&0&0&0&0&0&0&0\\0&0&1&0&0&0&0&0&0&0&0\\0
&0&0&1&0&0&0&0&0&0&0\\0&0&0&0&1&0&0&0&0&0&0\\0&0&0&0&0&\cos s&\sin
s&0&0&0&0\\0&0&0&0&0&-\sin s&\cos s&0&0&0&0\\0&0&0&0&0&0&0&\cos
s&\sin s&0&0\\0&0&0&0&0&0&0&-\sin s&\cos
s&0&0\\0&0&0&0&0&0&0&0&0&\cos s&-\sin s\\0&0&0&0&0&0&0&0&0&-\sin
s&\cos s\end{array}\right),
M^\varepsilon_2=\tiny\left(\begin{array}{ccccccccccc}
1&0&0&0&0&0&0&0&0&0&0\\0&1&0&0&0&0&0&0&0&0&0\\0&0&1&0&0&0&0&0&0&0&0\\0
&0&0&1&0&0&0&0&0&0&0\\0&0&-2a^2s&ks&1&0&0&0&0&0&0\\0&0&0&0&0&1&0&0&0&0&0\\
0&0&0&0&0&0&1&0&0&0&0\\0&0&0&0&0&0&-s&1&0&0&0\\0&0&0&0&0&s&0&0&1&0&0\\0&0&0&0&0&0&0&0&0&1&0\\
0&0&0&0&0&0&0&0&0&0&1\end{array}\right), \qquad \cdots\\
&\displaystyle \mbox{7 matrices} \quad\cdots\quad
M^\varepsilon_{11}=\tiny\left(\begin{array}{ccccccccccc} \cosh
as&0&0&0&0&0&0&0&0&\frac{1}{a}\sinh
as&0\\0&1&0&0&0&0&0&0&0&0&0\\0&0&\cosh
\sqrt{2}as&\frac{k}{\sqrt{2}a}(1-\cosh
\sqrt{2}as)&0&0&\frac{1}{\sqrt{2}a}\sinh \sqrt{2}as&0&0&0&0\\0
&0&0&1&0&0&0&0&0&0&0\\0&0&0&0&\cosh as&0&0&\sinh as&0&0&0\\0&0&0&0&0&1&0&0&0&0&0\\
0&0&\sqrt{2}a\sinh \sqrt{2}as&-\frac{k}{\sqrt{2}a}\sinh
\sqrt{2}as&0&0&\cosh
\sqrt{2}as&0&0&0&0\\0&0&0&0&\sinh s&0&0&\cosh s&0&0&0\\0&0&0&0&0&0&0&0&1&0&0\\a\sinh as&0&0&0&0&0&0&0&0&\cosh as&0\\
0&0&0&0&0&0&0&0&0&0&1\end{array}\right).
\end{eqnarray*}
by acting these matrices on a vector field $\textbf{v}$
alternatively we can  show that a one-dimensional optimal system
of ${\goth g}$ is given by
\begin{eqnarray*}\begin{array}{lll}
X_1=a_1\textbf{v}_1+a_2\textbf{v}_3+a_3\textbf{v}_4+a_4\textbf{v}_8&&
\hspace{-4cm}X_2=a_1\textbf{v}_1+a_2\textbf{v}_3+a_3\textbf{v}_4+a_4\textbf{v}_5-a_6\textbf{v}_9\\
X_3=a_1\textbf{v}_1+a_2\textbf{v}_2+a_3\textbf{v}_3+a_4\textbf{v}_4+a_6\textbf{v}_8&&
\hspace{-4cm}X_4=a_1\textbf{v}_1+a_2\textbf{v}_3+a_3\textbf{v}_4+a_4\textbf{v}_5-a_6\textbf{v}_8\\
X_5=a_1\textbf{v}_1+a_2\textbf{v}_3+a_3\textbf{v}_5+a_4\textbf{v}_6+a_6\textbf{v}_8&&
\hspace{-4cm}X_6=a_1\textbf{v}_1+\textbf{v}_2+a_2\textbf{v}_3+a_3\textbf{v}_4+a_4(\textbf{v}_6-\textbf{v}_{10})\\
X_7=a_1\textbf{v}_1+\textbf{v}_2+a_2\textbf{v}_3+a_3\textbf{v}_4+\textbf{v}_7-\textbf{v}_{11},&&
\hspace{-4cm}X_8=a_1\textbf{v}_1+a_2\textbf{v}_2+a_3\textbf{v}_3+a_4\textbf{v}_4+a_5\textbf{v}_5+a_6\textbf{v}_6,\\
X_9=a_1\textbf{v}_1+a_2\textbf{v}_2+\textbf{v}_3+\textbf{v}_4-(2a^2-k)\textbf{v}_5+\textbf{v}_7,&&
\hspace{-4cm}X_{10}=a_1\textbf{v}_1+a_2\textbf{v}_2+\textbf{v}_3+a_3\textbf{v}_4-(2a^2-k)\textbf{v}_5+\textbf{v}_6,\\
X_{11}=a_1\textbf{v}_1+a_2\textbf{v}_2+a_3\textbf{v}_3+a_4\textbf{v}_4+a_5\textbf{v}_5+a_6\textbf{v}_6+a_7\textbf{v}_{11},\\
X_{12}=a_1\textbf{v}_1+a_2\textbf{v}_2+a_3\textbf{v}_3+a_4\textbf{v}_4+a_5\textbf{v}_5+a_6\textbf{v}_6+a_7\textbf{v}_7+a_8\textbf{v}_{11},&&\\
X_{13}=a_1\textbf{v}_1+a_2\textbf{v}_2+a_3\textbf{v}_3+a_4\textbf{v}_4+a_5\textbf{v}_5+a_6\textbf{v}_6-a_7\textbf{v}_7-a_8\textbf{v}_9,&&\\
X_{14}=a_1\textbf{v}_1+a_2\textbf{v}_2+a_4\textbf{v}_3+a_4\textbf{v}_4+a_5\textbf{v}_5+a_6\textbf{v}_6+\textbf{v}_7-(2a+k)\textbf{v}_{11},&&\\
X_{15}=\frac{1}{2a^2-k}(\textbf{v}_1+\textbf{v}_8)+a_1\textbf{v}_2+a_2\textbf{v}_3+a_3\textbf{v}_4+a_5\textbf{v}_5+\textbf{v}_6+\textbf{v}_7,&&\\
X_{16}=a_1\textbf{v}_1+a_2\textbf{v}_2+\textbf{v}_3+\textbf{v}_4-(2a^2-k-1)\textbf{v}_5+a_3\textbf{v}_6+a_4\textbf{v}_7+a_5\textbf{v}_8+a_6\textbf{v}_9,&&
\end{array}\end{eqnarray*}
\section{Lie Algebra Structure}
$~~~~~$In this part, we determine the structure of symmetry Lie
algebra of the telegraph equation.\\
$\goth g$ has a \textit{Levi decomposition} in the form of ${\goth
g}={\goth r}\ltimes{\goth h}$, where ${\goth
r}=\langle\textbf{v}_2,\textbf{v}_3,\textbf{v}_4,\textbf{v}_5,\textbf{v}_6\rangle$
is the radical (the large solvable ideal) of $\goth g$ which is a
nilpotent nilradical of $\goth g$ and ${\goth
h}=\langle\textbf{v}_1,\textbf{v}_7,\textbf{v}_8,\textbf{v}_9,\textbf{v}_{10},\textbf{v}_{11}\rangle$
is a solvable and non-semisimple subalgebra of $\goth g$ with
centralizer $\langle\textbf{v}_4\rangle$ containing in the minimal
ideal
$\langle\textbf{v}_1,\textbf{v}_2,2a^2\textbf{v}_3-k\textbf{v}_4,\textbf{v}_5,...,\textbf{v}_{11}\rangle$.

Here we can find the quotient algebra generated from $\goth g$
such as
\begin{eqnarray}\label{eq:2}
{\goth g}_1={\goth g}/{\goth r},
\end{eqnarray}
with commutators table (\ref{table:2}), where
$\textbf{w}_i=\textbf{v}_i+{\goth r}$ for $i=1,...,11$ are members
of ${\goth g}_1$.

The (\ref{eq:2}) helps us to reduction differential equations. If
we want to integration an involutive distribution, the process
decomposes into two steps:
\begin{itemize}
\item integration of the involutive distribution with symmetry Lie
algebra ${\goth g}/{\goth r}$, and
\item integration on integral manifolds with symmetry algebra $\goth
r$.
\end{itemize}
$~~~~$First, applying this procedure to the radical $\goth r$ we
decompose the integration problem into two parts: the integration
of the distribution with semisimple algebra ${\goth g}/{\goth r}$,
then the integration of the restriction of distribution to the
integral manifold with the solvable symmetry algebra $\goth r$.\\
$~~~~$The last step can be performed by quadratures. Moreover,
every semisimple Lie algebra ${\goth g}/{\goth r}$ is a direct sum
of simple ones which are ideal in ${\goth g}/{\goth r}$. Thus, the
Lie-Bianchi theorem reduces the integration problem ti involutive
distributions equipped with simple algebras of symmetries.

$~~~~$Both $\goth g$ and ${\goth g}_1$ are non-solvable, because
if ${\goth
g}^{(1)}=\langle\textbf{v}_i,[\textbf{v}_i,\textbf{v}_j]\rangle=[\goth
g, \goth g]$, and ${\goth
g}_1^{(1)}=\langle\textbf{w}_i,[\textbf{w}_i,\textbf{w}_j]\rangle=[{\goth
g}_1,{\goth g}_1]$, be the derived subalgebra of $\goth g$ and
${\goth g}_1$ we have
\begin{eqnarray*}
{\goth g}^{(1)}=[{\goth g},{\goth g}] =\langle
\textbf{v}_1,\textbf{v}_2,2a^2\textbf{v}_3-k\textbf{v}_4,\textbf{v}_5,\textbf{v}_6,\textbf{v}_7,\textbf{v}_8,\textbf{v}_9,
\textbf{v}_{10},\textbf{v}_{11}\rangle=[{\goth g}^{(1)},{\goth
g}^{(1)}] ={\goth g}^{(2)},
\end{eqnarray*}
and
\begin{eqnarray*}
{\goth g}_1^{(1)}=[{\goth g}_1,{\goth g}_1]=\langle
\textbf{w}_1,\textbf{w}_2,\textbf{w}_3,\textbf{w}_4,\textbf{w}_5,\textbf{w}_6\rangle
={\goth g}_1.
\end{eqnarray*}
Thus, we have a chain of ideals ${\goth g}\supset{\goth
g}^{(1)}={\goth g}^{(2)}\neq 0$, ${\goth g}_1\supset{\goth
g}_1^{(1)}={\goth g}_1\neq 0$, which sows the non-solvability of
$\goth g$ and ${\goth g}_1$.
\begin{table}
\caption{Commutation relations of $\goth g$ }\label{table:2}
\vspace{-0.3cm}\begin{eqnarray*}\hspace{-0.75cm}\begin{array}{ccccccc}
\hline
  [\,,\,]        &\hspace{1cm}\textbf{w}_1       &\hspace{0.5cm}\textbf{w}_2    &\hspace{0.5cm}\textbf{w}_3    &\hspace{0.5cm}\textbf{w}_4    &\hspace{0.5cm}\textbf{w}_5       &\hspace{0.5cm}\textbf{w}_6  \\ \hline
  \textbf{w}_1   &\hspace{1cm} 0                 &\hspace{0.5cm} 0              &\hspace{0.5cm}-\textbf{w}_4   &\hspace{0.6cm}\textbf{w}_3    &\hspace{0.5cm}-\textbf{w}_6      &\hspace{0.5cm}\textbf{w}_5 \\
  \textbf{w}_2   &\hspace{1cm} 0                 &\hspace{0.5cm} 0              &\hspace{0.7cm}\textbf{w}_6    &\hspace{0.5cm}-\textbf{w}_5   &\hspace{0.5cm}-\textbf{w}_4      &\hspace{0.5cm}\textbf{w}_3 \\
  \textbf{w}_3   &\hspace{0.9cm} \textbf{w}_4    &\hspace{0.3cm} -\textbf{w}_6  &\hspace{0.7cm}0               &\hspace{0.5cm}-\textbf{w}_1   &\hspace{0.7cm}0                  &\hspace{0.5cm}\textbf{w}_2 \\
  \textbf{w}_4   &\hspace{0.8cm} -\textbf{w}_3   &\hspace{0.4cm} \textbf{w}_5   &\hspace{0.7cm}\textbf{w}_1    &\hspace{0.7cm}0               &\hspace{0.5cm}-\textbf{w}_2      &\hspace{0.5cm}0 \\
  \textbf{w}_5   &\hspace{0.9cm} \textbf{w}_6    &\hspace{0.4cm}\textbf{w}_4    &\hspace{0.7cm}0               &\hspace{0.7cm}\textbf{w}_2    &\hspace{0.7cm}0                  &\hspace{0.4cm}a^2\textbf{w}_1 \\
  \textbf{w}_6   &\hspace{0.8cm} -\textbf{w}_5   &\hspace{0.3cm} -\textbf{w}_3  &\hspace{0.5cm}-\textbf{w}_2   &\hspace{0.7cm}0               &\hspace{0.5cm}-a^2\textbf{w}_1   &\hspace{0.5cm}0 \\
  \hline\end{array}\end{eqnarray*}\end{table}
\section{Conclusion}
In this article group classification of telegraph equation and the
algebraic structure of the symmetry group is considered.
Classification of one-dimensional subalgebra is determined by
constructing one-dimensional optimal system. The structure of Lie
algebra symmetries is analyzed.

%\bibliographystyle{plain}        % Include this if you use bibtex
%\bibliography{autosam}           % and a bib file to produce the
                                 % bibliography (preferred). The
                                 % correct style is generated by
                                 % Elsevier at the time of printing.

%\begin{thebibliography}{99}     % Otherwise use the
                                 % thebibliography environment.
                                 % Insert the full references here.
                                 % See a recent issue of Automatica
                                 % for the style.
%  \bibitem[Heritage, 1992]{Heritage:92}
%     (1992) {\it The American Heritage.
%     Dictionary of the American Language.}
%     Houghton Mifflin Company.
%  \bibitem[Able, 1956]{Abl:56}
%     B.~C.~Able (1956). Nucleic acid content of macroscope.
%     {\it Nature 2}, 7--9.
%  \bibitem[Able {\em et al.}, 1954]{AbTaRu:54}
%     B.~C. Able, R.~A. Tagg, and M.~Rush (1954).
%     Enzyme-catalyzed cellular transanimations.
%     In A.~F.~Round, editor,
%     {\it Advances in Enzymology Vol. 2} (125--247).
%     New York, Academic Press.
%  \bibitem[R.~Keohane, 1958]{Keo:58}
%     R.~Keohane (1958).
%     {\it Power and Interdependence:
%     World Politics in Transition.}
%     Boston, Little, Brown \& Co.
%  \bibitem[Powers, 1985]{Pow:85}
%     T.~Powers (1985).
%     Is there a way out?
%     {\it Harpers, June 1985}, 35--47.

%\end{thebibliography}

%\appendix
%\section{A summary of Latin grammar}    % Each appendix must have a short title.
%\section{Some Latin vocabulary}         % Sections and subsections are supported
                                        % in the appendices.
\end{document}